\documentclass[12pt]{amsart}

\usepackage{amsthm}

\usepackage{graphicx}
\usepackage{epsfig}
\usepackage{epstopdf}

\usepackage{latexsym}
\usepackage{amsmath}
\usepackage{amssymb}
\usepackage{amsfonts}

\theoremstyle{definition}

\newtheorem{dfn}{Definition}[section]

\newtheorem{example}[dfn]{Example}
\newtheorem{rem}[dfn]{Remark}
\newtheorem{remark}[dfn]{Remark}
\theoremstyle{plain}

\newtheorem{problem}[dfn]{Problem}

\newtheorem{thm}[dfn]{Theorem}

\newtheorem{lem}[dfn]{Lemma}
\newtheorem{lemma}[dfn]{Lemma}

\newtheorem{prop}[dfn]{Proposition}
\newtheorem{proposition}[dfn]{Proposition}

\newtheorem{assumption}[dfn]{Assumption}
\newtheorem{cor}[dfn]{Corollary}

\newtheorem{question}[dfn]{Question}

\def\lra{\longrightarrow}
\def\proof{\par\medskip\noindent{\it Proof. }}

\def\k{{\mathbf k}}
\def\R{{\mathbb R}}

\def\C{{\mathbb C}}
\def\Z{{\mathbb Z}}
\def\P{{\mathbb P}}

\def\Q{{\mathbb Q}}

\def\N{{\mathbb N}}

\def\eps{\epsilon}
\def\al{\alpha}
\def\be{\beta}
\def\ga{\gamma}
\def\Ga{\Gamma}

\def\De{\Delta}

\def\Si{\Sigma}

\def\La{\Lambda}
\def\Om{\Omega}
\def\om{\omega}

\def\acts{\curvearrowright}
\def\D{\partial}

\def\embed{\hookrightarrow}
\def\8{\infty}
\def\<{\langle}
\def\>{\rangle}

\def\ol{\overline}

\renewcommand{\hom}{\operatorname{Hom}}
\newcommand{\h}{\operatorname{H}}

\def\t{\tilde}

\usepackage{epsfig}

\long\def\comment#1\endcomment{}

\input diagrams.tex

\begin{document}

\title{	On representation varieties of 3-manifold groups}
\author{Michael Kapovich, John J. Millson}
\date{\today}

\begin{abstract}
We prove universality theorems (``Murphy's Laws") for representation varieties of fundamental groups of closed 3-dimensional manifolds. We show that germs of $SL(2)$-representation schemes of such groups are essentially the same as germs of schemes of finite type over $\Q$. 
\end{abstract}

\maketitle

\section{Introduction}

In this paper we will prove that there are no ``local'' restrictions on geometry of representation schemes 
of 3-manifold groups to $PO(3)$ and $SL(2)$. Note that both groups $H=PO(3)$ and $H=SL(2)$ are affine algebraic group schemes defined over $\Q$, thus, for every finitely-generated group $\Gamma$, the representation schemes 
$$
\hom(\Ga, H)
$$
and character schemes 
$$
X(\Ga, H)= \hom(\Ga, H)//H
$$
are affine algebraic schemes over $\Q$. Our goal is to show that, {\em to some extent}, these are the only restrictions on local geometry of the representation and character schemes of fundamental groups of closed 3-manifolds. The universality theorem we thus obtain is one in many universality theorems about moduli spaces of geometric objects, see \cite{Mnev}, \cite{RG}, \cite{KM}, \cite{KM2}, \cite{KM3}, \cite{V}, \cite{P}, 
\cite{Rapinchuk}.

Below is the precise formulation of our universality theorem. In what follows we will be using the notation 
$G=PO(3), \tilde G= Spin(3)$.

\begin{thm}\label{main}
Let $X\subset \C^N$ be an affine algebraic scheme over $\Q$ and let $x\in X$ be a rational   
point. Then there exist:

1. An open subscheme $X'\subset X$ containing $x$.

2. A closed 3-dimensional  manifold $M$ with the fundamental group $\pi$. 

3. A representation $\rho_0: \pi\to PO(3, \R)$, such that the image of $\rho_0$ is dense in $PO(3,\R)$. 

4. An open $G$-invariant subscheme $R'\subset \hom(\pi, G)$ containing $\rho_0$ 
and a closed subscheme $R'_c\subset R'$ which is a cross-section for the action 
$$
G\times R'\to R'. 
$$

5. An isomorphism of schemes over $\Q$: 
$$
f: R'\to G^k \times X', \quad f(\rho_0)=(1, x),$$ 
for some $k$. 
\end{thm}

\begin{rem}
One can show that the same theorem holds for a homomorphism $\rho_0$ whose image is a finite group with trivial centralizer in $PO(3, \R)$. 
\end{rem}

Since the groups $PSL(2, \C)$ and $PO(3, \C)=G(\C)$ are isomorphic, and 
$$\t G(\C)=Spin(3, \C)\cong SL(2, \C),$$
``universality'' for $PO(3)$--representations  leads to the one for $SL(2)$--represen\-tations:

\begin{cor}\label{cor:main}
Let $X\subset \C^N$ be an affine algebraic scheme over $\Q$ and $x\in X$ be a rational 
point. Then there exist:

1. An open subscheme $X'\subset X$ containing $x$.

2. A closed 3-dimensional  manifold $M$ with the fundamental group $\pi$. 

3. A representation $\tilde\rho_0: \pi\to SU(2)< SL(2, \C)$, such that the image of $\t\rho_0$ is dense in $SU(2)$. 

4. An open $G$-invariant subscheme $\widehat{S}\subset \hom(\pi, SL(2))$ containing $\t\rho_0$,  
such that every complex point of $\widehat{S}$ is a Zariski dense representation to $SL(2,\C)$. 

5. A regular  \'etale covering of schemes over $\C$: 
$$
\t{f}: \widehat{S} \to \t G ^{k} \times  X', \quad \t{f}(\tilde\rho_0)=(1, x),
$$
with the deck-transformation group isomorphic to $\Z_2^r$, for some $k$ and $r$.  
\end{cor} 

In particular, $\t{f}$ yields an isomorphism of the analytic germs 
$$
(\hom(\pi, SL(2)), \t\rho_0)\to (\C^{3k+3} \times X', 0\times x),
$$
for some $k\ge 1$. Thus, if the scheme $X'$ is non-reduced at $x$, so is $\hom(\pi, SL(2))$.

\begin{rem}
Despite of our efforts, we were unable to replace an  \'etale covering with an isomorphism in Corollary \ref{cor:main}. This is strangely reminiscent of the finite abelian coverings appearing in our universality theorem for planar linkages, \cite{KM3}. Note that a relation between universality theorems for projective arrangements  and spherical linkages was established in \cite{KM2}, where a finite abelian covering appeared for essentially the same reason as in the present paper.  
\end{rem}

We will see that the action of  $SL(2)$ on $\t R'$ factors through  the group 
$PSL(2)$, which admits a cross-section.  In particular, we obtain

\begin{cor}\label{cor:germs}
There exists an isomorphism of analytic germs
$$
 (X(\pi, SL(2)), [\t\rho_0]) \to (\C^{3k}\times X', 0\times x).
$$
\end{cor}

\begin{example}
Pick a natural number $\ell$. Then there exists 
a closed 3-dimensional manifold $M$, an integer $n$ and a representation $\rho: \pi_1(M)\to SU(2)$ with dense image, so that the completed local ring of the germ  
$$
X(\pi_1(M), SL(2)), [\rho])
$$ 
is isomorphic to the completion of the ring
$$
\C[t, t_1,\ldots, t_{3k}]/(t^\ell).  
$$
\end{example}

This shows that the representation and character schemes of 3-manifold groups can be nonreduced (at points of Zariski density), which is why we refrain from referring to these schemes as ``varieties,'' as it is commonly done in the literature.

\begin{remark}
Recently Igor Rapinchuk proved a universality theorem for character schemes of groups $\Gamma$ satisfying Kazhdan's Property T: It involves representations of $\Gamma$'s into $SL(n,\C)$, \cite[Theorem 3]{Rapinchuk}. Unlike \cite{KM} and this paper, Rapinchuk's theorem applies to the entire character variety $X^{red}(\Ga, SL(n,\C))$ minus the trivial representation (which is an isolated point). In Rapinchuk's theorem, the number $n$ (and the group $\Ga$) depend on the given affine variety $X$ over $\Q$. 
\end{remark}

\newtheorem*{ack}{Acknowledgments}

 \begin{ack} 
{\em Partial financial support  to the first author was provided by the NSF grant DMS-12-05312 
and to the second author by the NSF grant DMS-12-06999. We are grateful to the referee for useful remarks.} 
\end{ack}

\section{Preliminaries}

\subsection{Representation and character schemes}

We will say that a subscheme $Y\subset X$ is {\em clopen} if it is both closed and open. We will use the topologist's notation:
$$
\Z_m:= \Z/m \Z,
$$
is the cyclic group of order $m$. 

Let $G$ be an algebraic group scheme over a field $\k$ 
of characteristic zero (this will be the default assumption through the rest of the paper) 
with the Lie algebra ${\mathfrak g}$. Let  $\Ga$ be a finitely-presented group with the presentation
$$
\<s_1,\ldots, s_p | r_1=1, \ldots , r_q=1\>.
$$
(In fact, one needs $\Ga$ only to be finitely-generated, but all finitely-generated 
groups in this paper will be also finitely-presented.) 
Every word $w$ in the generators  $s_i, s_i^{-1}, i=1,\ldots,p$,  defines a morphism 
$$
w: G^p\to G,
$$
obtained by substituting elements $g_1^{\pm 1},\ldots,g_p^{\pm 1}\in G$ 
in the word $w$ for the letters 
$s_1^{\pm 1},\ldots,s_p^{\pm 1}$. We then obtain the {\em representation scheme}
$$
\hom(\Ga, G)=\{ (g_1,\ldots, g_p)\in G^p: r_j(g_1,\ldots, g_p)=1, j=1,\ldots,q\}, 
$$
as every homomorphism $\Ga\to G$ is determined by its values on the generators of $\Ga$. We will, thus, think of points of this scheme as homomorphisms $\rho: \Ga\to G$. 
The representation scheme is known to be independent of the presentation of the group $\Ga$. We refer the reader to \cite{LM, Si} for more details. We also refer the reader to \cite{S, Si} for detailed discussion of character varieties/schemes and survey of their applications to 3-dimensional topology. 

We will frequently use the following two facts about representation schemes, see e.g \cite{Si}:

1. $\hom(\Ga_1\star \ldots \star \Ga_k, G)\cong \prod_{i=1}^k \hom(\Ga_i, G)$. 

\smallskip 
2. For each $\rho\in \hom(\Ga, G(\k))$ satisfying $\h^1(\Ga, {\mathfrak g}_{Ad\rho})=0$, 
$$
\rho\in \hom(\Ga, G(\k))$$ 
is a smooth point of the scheme $\hom(\Ga, G)$. The $G$-orbit through $\rho$ is open in 
$\hom(\Ga, G)$. 

\medskip
In what follows we will use the simplified  
notation $\h^q(\Ga, {Ad\rho})$ instead of $\h^q(\Ga, {\mathfrak g}_{Ad\rho})$.

\medskip 
We assume from now on that the group $G$ is affine; in particular, $\hom(\Ga, G)$ is also 
an affine scheme. The group $G$ acts naturally on this scheme:
$$
(g, \rho)\mapsto \rho^g, \quad \rho^g(\ga)= g \rho(\ga) g^{-1}. 
$$
Assuming, in addition, that $G$ is reductive, we obtain the GIT quotient
$$
X(\Ga, G)= \hom(\Ga, G)//G,
$$
which is a scheme of finite type known as the {\em character scheme} (or, more commonly, as the {\em character variety}). However, as we will see, both representation and character schemes can be  nonreduced, so we will avoid the traditional {\em representation/character variety} terminology. 

We will use the notation 
$$
\hom^{red}(\Ga, G), \quad X^{red}(\Ga, G)
$$
to denote the varieties which are the reductions of the schemes 
$$
\hom(\Ga, G), \quad X(\Ga, G). 
$$

\medskip
Recall that for every $\rho\in \hom(\Ga, G)$, the vector space of cocycles 
$$
Z^1(\Ga, Ad \rho)  
$$   
is isomorphic to the Zariski tangent space $T_\rho \hom(\Ga, G)$ and this isomorphism carries the 
subspace of coboundaries $B^1(\Ga, Ad\rho)$ to the tangent space of the $G$-orbit through $\rho$. 
Note, however, that  $\h^1(\Ga, Ad\rho)$ is {\em not} always isomorphic to the Zariski tangent space 
of $[\rho]\in X(\Ga, G)$, see \cite[\S 6]{H} as well as \cite{Si}. 

\medskip
Suppose now that the group $\Phi$ is finite.  
Then for every $\rho\in \hom(\Phi, G)$, 
$$
\h^1(\Phi, Ad \rho)=0. 
$$   
(Furthermore, $\h^i(\Phi, Ad \rho)=0, i\ge 1$.) In particular, the $G$-orbit of $\rho$ 
is a clopen (closed and open) subscheme 
$$
\hom_\rho(\Phi, G)\subset \hom(\Phi, G). 
$$
This subcheme is isomorphic to the quotient $G/\zeta_G(\rho(\Phi))$, where $\zeta_G(H)$ denotes  the centralizer of the subgroup $H$ in $G$. (Note that if $\zeta_G(\rho(\Phi))$ equals the center of $G$, then the point $[\rho]\in X(\Phi, G)$ is a reduced isolated point in the character scheme and the entire character scheme is smooth.) We obtain:

\begin{lem}\label{lem:L1}
For every finite group $\Phi$ and connected affine group $G$, the scheme $\hom(\Phi, G)$ is smooth and each irreducible component of $\hom(\Phi, G)$ is $G$-homogeneous. These irreducible components are the open subschemes  
$\hom_\rho(\Phi, G)$. If the representation $\rho$ is trivial, then $\hom_\rho(\Phi, G)$ is a single point.   
\end{lem}

The following lemma is also immediate: 

\begin{lemma}\label{lem:L2}
Let $\phi: \Ga\to \Ga'$ be a group homomorphism. Then the pull-back map 
$\phi^*(\rho)=\rho\circ \phi$ is a morphism of schemes
$$
\hom(\Ga', G)\to \hom(\Ga, G). 
$$
\end{lemma}

\begin{lemma}\label{lem:L3}
Let $\Ga$ be a finitely-presented group and let $\Theta\subset \Gamma$ be a finite subset with the quotient group 
$$
\Gamma':= \Gamma/\<\<\Theta\>\>.$$
Let $\phi: \Ga\to \Ga'$ denote the projection homomorphism. Then the pull-back morphism 
$$
\phi^*: \hom(\Ga', G)\to \hom_\Theta(\Ga, G) 
$$
is an isomorphism, where 
$$
\hom_\Theta(\Ga, G) \subset \hom(\Ga, G)
$$
is the closed subscheme defined by
$$
\hom_\Theta(\Ga, G)=\{\rho\in \hom(\Ga, G): \rho(\theta)=1, \forall \theta\in \Theta\}. 
$$
\end{lemma}
\proof Given a finite presentation $P$ of $\Ga$ let $P'$ be the presentation of $\Ga'$ obtained from $P$ by 
adding words representing elements of $\Theta$ as the relators. Then the assertion follows immediately from the definition of the representation scheme of a group using a group presentation. 

\begin{cor}\label{cor:C3}
Suppose that every element $\theta\in \Theta$ has finite order. Then
the isomorphism $\phi_*:  \hom(\Ga', G)\to \hom_\Theta(\Ga, G)$ 
 sends $\hom(\Ga', G)$ to the open  subscheme $\hom_\Theta(\Ga, G)\subset 
 \hom(\Ga, G)$. \end{cor}
\proof Consider an element $\theta\in \Theta$ and the trivial representation $\rho_{0,\theta}: \<\theta\>\to G$. 
By Lemma \ref{lem:L1}, the singleton  
$$
\{\rho_{0,\theta}\}=\hom_{\rho_0,\theta}(\<\theta\>, G) \subset \hom(\<\theta\>, G)
$$
is a reduced isolated point in the scheme $\hom(\<\theta\>, G)$. In particular, this point is open in 
$\hom(\<\theta\>, G)$. We have the pull-back morphism
$$
\phi^*_{\theta}: \hom(\Ga, G)\to \hom(\<\theta\>, G),
$$
induced by the inclusion homomorphism $\phi_{\theta}: \<\theta\>\embed \Ga$. Therefore,
$$
(\phi^*_{\theta})^{-1}\left( \hom_{\rho_0,\theta}(\<\theta\>, G)\right) \subset \hom(\Ga, G)
$$
is an open subscheme. Furthermore, by the definition of $\hom_\Theta(\Ga, G)$, 
$$
\hom_\Theta(\Ga, G) =\bigcap_{\theta\in \Theta} (\phi^*_{\theta})^{-1}\left( \hom_{\rho_0,\theta}(\<\theta\>, G)\right).
$$
(A homomorphism belongs to $\hom_\Theta(\Ga, G)$ if and only if it sends each $\theta\in \Theta$ to $1\in G$.) Therefore, $\hom_\Theta(\Ga, G)$ is also open in $\hom(\Ga, G)$. \qed

\subsection{Coxeter groups}\label{sec:Coxeter}

We refer the reader to \cite{Humphreys} for the basics of Coxeter groups. 

Let $\Delta$ be a finite simplicial graph with the vertex and edge sets denoted $V=V(\Delta)$ and $E=E(\Delta)$ respectively. We will use the notation $e=[v,w]$ for the edge of $\Delta$ connecting $v$ and $w$, if it exists.  
We assume also that we are given a function 
$$
m: E\to \N, \quad m(e)\ge 2, \forall e\in E, 
$$ 
{\em labeling the edges} of $\De$. We will call the pair $(\Delta, m)$ a {\em labeled graph} or a {\em Coxeter graph}. 
Given a labeled graph, we define the associated {\em Coxeter group} $\Ga=\Ga_\Delta$ by the presentation
$$
\< g_v, v\in V| \forall v, w\in V, g_v^2=1, \underbrace{g_v g_w \cdots}_{m \hbox{~terms}}=
\underbrace{g_w g_v \cdots}_{m \hbox{~terms}}, ~~e=[v,w], m=m(e), e\in E\>.  
$$
Alternatively, one can describe the relators of this group as $g_v^2=1$ and
$$
(g_vg_w)^m=1
$$
where $m=m(e), e=[v,w]$. 

\begin{rem}
Note that the notation we use here is different from the one in the Lie theory, 
where two generators commute whenever 
the corresponding vertices are not connected by an edge. In our notation, every such pair of 
elements of $\Ga$ generates an infinite dihedral subgroup of $\Ga$. 
\end{rem}

We also define the {\em canonical central extension}
\begin{equation}\label{eq:cen}
1\to \Z_2\to \t\Ga \stackrel{\eta}{\to} \Ga \to 1
\end{equation}
of the group $\Ga$, with  the {\em extended Coxeter group} $\t\Ga=\t\Ga_\Delta$ given by the presentation
$$
\< z, g_v, v\in V| z^2=1, \forall v\in V, [g_v, z]=1, g_v^2=z,$$
$$
\underbrace{g_v g_w \cdots}_{m \hbox{~terms}}=z^{m+1} \underbrace{g_w g_v \cdots}_{m \hbox{~terms}}, ~~e=[v,w], m=m(e), e\in E\>. 
$$
The number $r=|V|$ (the cardinality of $V$) is called the {\em rank} of $\Ga$ and $\t\Ga$. We will refer to the generator $z$ of the group $\Ga$ as {\em the central element} of $\Ga$, even though, the center of $\t\Ga$ might be  larger than $\Z_2$: This happens precisely when $\De$ consists of a single vertex.  

\medskip
A subgraph $\Si\subset \Delta$ is called {\em full} if for every vertices $v, w\in \Si$, the edge 
$[v,w]$ in $\Delta$ also belongs to $\Sigma$. Every subgraph $\Si\subset \Delta$ inherits labels from $\Delta$.  For the new labeled graph (which we still denote $\Sigma$), we have the natural homomorphism,
$$
\iota_{\Si}: \Ga_\Si\to \Ga_\Delta
$$
sending each generator $g_v\in \Ga_\Si$, $v\in V(\Si)$, to the generator of $\Ga_\Delta$ with the same name. 
 It is immediate that the homomorphism $\iota_\Si$ lifts to a homomorphism
$$
\tilde \iota_\Si: \t\Ga_\Si\to \t\Ga_\Delta
$$
sending each $g_v$ to itself ($v\in V(\Si)$) and the central element $z\in \t\Ga_\Si$ to the central element 
$z\in \t\Ga_\Delta$. We will use this construction in two spacial cases:

a. $\Si:=\Delta_{\emptyset}$ is the subgraph which has the same vertex set as $\Ga$, but empty edge set. 
Then 
$$
\Ga_\Si\cong F_r, \quad \t\Ga_\Si\cong F_r \times \Z_2. 
$$

b. $\Si\subset \Delta$ is a full subgraph. In this case,   
the homomorphism $\iota_\Si$  is injective, see e.g. \cite{Humphreys},  page 113; it follows that the homomorphism $\tilde \iota_\Si$ is injective as well.

\medskip 
For full subgraphs $\Si\subset \Delta$, 
the subgroups $\iota_\Si(\Ga_\Si)< \Ga_\Delta$ and $\tilde\iota_\Si(\tilde\Ga_\Si)< \tilde\Ga_\Delta$
are called {\em parabolic} subgroups of $\Ga_\Delta$ and $\t\Ga_\Delta$ respectively. 
We say that a parabolic subgroup of $\Ga_\Delta, \t\Ga_\Delta$ is {\em elementary}, if it is a finite 
parabolic subgroup of rank $\le 2$. The latter requirement simply means that $\Si$ consists of at most two vertices; the finiteness condition means that if $\Si$ consists of two vertices, then these vertices are connected by an edge. We will refer to such subgraphs $\Si$ {\em elementary} as well. 

\begin{example}\label{e1}
1. Suppose that $\Delta$ consists of a single edge $e$ labelled $2$. 
Then $\Ga_\Delta\cong \Z_2\times \Z_2$ and 
$$
\t\Ga_\Delta \cong Q_8,
$$ 
the finite quaternion group. 

2. If the edge $e$ is labeled $4$ then $\Ga_\Delta$ is the dihedral group $I_2(4)$ of order $8$; it admits an epimorphism 
$$
\Ga_\Delta\to \Z_2\times \Z_2,
$$
whose kernel is the center of $\Ga_\Delta$, which is generated by the involution $g_v g_w g_v g_w$.
\end{example}

\section{Representations of Coxeter groups and extended Coxeter groups}

In this section we prove some basic facts about representations of Coxeter and 
extended Coxeter groups to $PSL(2,\C)$ and $SL(2,\C)$ respectively.

\subsection{Representations of elementary Coxeter groups}

Recall that the quotient map
$$
p: \t G(\C)=SL(2,\C)\to G(\C)= PSL(2, \C)= SL(2,\C)/\{\pm 1\},
$$
is a 2-fold covering. The extended Coxeter groups appear naturally in the context of lifting homomorphisms of Coxeter groups from $PSL(2,\C)$ to $SL(2,\C)$. 

Consider the labelled graph $\Delta$,  
consisting of two vertices $v, w$ and the edge $[v,w]$ labelled $n\ge 2$. 
The corresponding Coxeter group $\Ga_\Delta$ is a finite dihedral group, usually denoted $I_2(n)$. This group is isomorphic to the subgroup of the group of symmetries of a regular planar $2n$-gon, acting simply transitively on the set of edges of this polygon. Hence, this group embeds canonically (up to conjugation) into $O(2)$ and, thus, into $PO(3, \R)\cong SO(3, \R) < PSL(2,\C)$. If $n$ is odd, then such group of symmetries will lift isomorphically to a subgroup of $SU(2)< SL(2, \C)$. In contrast, we will be interested (only) in the case when $n$ is even; in fact, we will be using only the labels $n=2$ and $n=4$ in this paper.

 Below we will prove several lemmas about faithful representations of elementary Coxeter (and extended Coxeter) groups into $G(\C)$ (and $\t G(\C)$). 

\medskip

\begin{lem}\label{lem:Z2timesZ2}
1. There exists, unique up to conjugation, faithful representation \newline $\rho: \Z_2\times \Z_2\to G(\C)$. 

2. There are no faithful representations $\Z_2\times \Z_2\to \t G(\C)$.
\end{lem}
\proof 1. Since the image of $\rho$  is finite, its conjugate to a subgroup of $SO(3,\R)< G(\C)$. 
Part 1 now follows from the fact that the group $SO(3,\R)$ acts transitively on the set of pairs of orthogonal 1-dimensional subspaces of $\R^3$ (these subspaces, in our setting, are fixed lines of the images of the direct factors of $\Z_2\times \Z_2$ under $\rho$). 

2. Part 2 follows from the fact that any involution $A\in SL(2,\C)$ has both eigenvalues equal to $-1$, i.e., $A$ equals $-1\in SL(2,\C)$. \qed 

\medskip
The next lemma and proposition generalize Lemma \ref{lem:Z2timesZ2} to representations of the dihedral group $I_2(4)$.

\begin{lem}\label{lem:I_2(4)}
All injective representations $\rho: \Ga=I_2(4)\to G(\C)$ are conjugate to each other.  
\end{lem}
\proof Since the group $\Ga$ is finite, its image in $G(\C)$ lies in a conjugate of the maximal compact subgroup $SO(3, \R)< G(\C)$. Thus, we will assume that $\rho(\Ga)$ is contained in $SO(3,\R)$. Since the product of the generating involutions $\rho(g_v), \rho(g_w)$ 
of $\rho(\Ga)$ has order $4$, the fixed lines of $\rho(g_v), \rho(g_w)$ meet at the angle $\frac{\pi}{4}$ in $\R^3$. Now, the assertion follows from the fact that $SO(3, \R)$ acts transitively on the set of 1-dimensional subspaces in $\R^3$ meeting at the given angle. \qed

\begin{proposition}\label{prop:P0}
Consider the dihedral group $I_2(2m)=\Ga_\Delta$, and its isomorphism 
$$
\rho: \Ga_\Delta \to \Ga< G(\C). 
$$ 
Then:

1.  For every choice of matrices $\t g_u\in \t G(\C)$ projecting to the generators $\rho(g_u)\in \Ga < G(\C)$,  the map 
$$
g_u\to \t g_u, \quad u\in \{v, w\}, 
$$
extends to a monomorphism $\tilde\rho: \t \Ga_\Delta \to \t G(\C)$. 

2. The centralizer of the group $\t\rho(\t\Ga_\Delta)$ in $\t G(\C)$ equals the center of $\t G(\C)$.  
\end{proposition}
\proof The proof amounts to elementary linear algebra, we include the details for the sake of completeness. For the notational convenience we will identify the isomorphic groups $\Ga$ and $\Ga_\Delta$. After conjugating the subgroup $\Ga$ in $G(\C)$, we can (and will) assume that $\Ga$ lies in the subgroup 
$SO(3, \R)< G(\C)$. The orthogonal subgroup is covered by the unitary subgroup $SU(2)<SL(2,\C)$. We let 
$Z(SU(2))\cong \Z_2$ denote the center of $SU(2)$; this center consists of the matrices $\pm 1$.

We begin with several trivial observations. Since $\rho$ is injective, the involutions $g_v, g_w$ are distinct rotations in $SO(3,\R)$. In particular, their fixed-point sets in $\C P^1$ are pairwise disjoint. Suppose that the elements $\t g_u, \t g_v\in SU(2)$ project to $g_v, g_w$ respectively. Since the kernel of the covering 
$\t G\to G$ is isomorphic to $\Z_2$, the unitary transformations $\t g_v, \t g_w\in SU(2)$ have order at most $4$: 
$$
\t g_u^2\in Z(SU(2)), u\in \{v, w\}. 
$$ 
Note that the only involution in $SU(2)$ is the matrix $-1$. Since $\t g_u$ projects nontrivially to $SO(3)$, this matrix cannot be an involution. It follows that 
$$
\t g_u^2=-1\in SU(2), \quad u\in \{v, w\}. 
$$
The eigenvalues of the matrices $\t g_v, \t g_w$ have to be roots of unity of the order $4$, which implies that the spectrum of each matrix $\t g_u, u\in \{v, w\}$, equals $\{i, -i\}$. 

We next claim that the eigenspaces of unitary transformations $\t g_v, \t g_w$ are pairwise distinct. If not, then these matrices would be simultaneously diagonalizable, which would imply that their projections to $PSL(2,\C)$ are equal. (Two involutions in $SO(3,\R)$ which have same fixed point sets have to be the same.) 

Suppose now that $A\in SL(2,\C)$ is a matrix centralizing the subgroup $\<\t g_v, \t g_w\>$ 
generated by $\t g_v, \t g_w$. We claim that $A$ is a scalar matrix, i.e., an element of the center of $SL(2,\C)$. Indeed, since $A$ commutes with both $\t g_v, \t g_w$, it has to preserve the eigenspaces of each matrix $\t g_v, \t g_w$. (Here we are using the fact that the eigenvalues of $\t g_u$ are distinct, $u\in \{v, w\}$.) However, a nonscalar matrix in $SL(2,\C)$ cannot have three distinct eigenlines. Therefore, $A$ is a scalar matrix. This implies the second claim of the lemma. 

The generators $g_v, g_w$ satisfy 
$$
t=(g_v g_w)^m= (g_w g_v)^m,
$$
where $t$ is an order $2$ element, which belongs to the center of $\Ga$. (In the geometric realization of 
$\Ga_\Si$ as a group of symmetries of a regular $2n$-gon, the element $t$ corresponds to the order $2$ rotation, the central symmetry of the polygon.) 

If we had the relation
$$
(\t g_v \t g_w)^{m}= (\t g_w \t g_w)^{m},
$$
it would result in the monomorphism
$$
\al: \Ga\to SU(2), \quad \al(g_u)=\tilde g_u, u\in \{v, w\}, 
$$
lifting the  embedding $\rho: \Ga\embed SO(3, \R)$. The image of the center $Z(\Ga)$ of $\Ga$ would then be in the center of $\al(\Ga)$, hence, as we noted above, in the center of $SU(2)$. 
Then, the composition $\rho=p\circ \al$ would send $Z(\Ga)$ to $1$, which is a contradiction. 

This leaves us with the only possibility  
$$
(\t g_v \t g_w)^{m}= - (\t g_w \t g_v)^{m}. 
$$
To conclude, the map 
$$
g_v\mapsto \t g_v, \quad z\mapsto -1\in SL(2,\C),
$$
extends to an homomorphism $\t\Ga_\Delta\to p^{-1}(\Ga)$, sending the central element $z\in \t\Ga_\De$ to the matrix $-1\in SL(2,\C)$. Injectivity of this homomorphism follows from injectivity of the representation $\Ga\to PSL(2,\C)$. \qed

\subsection{Representations faithful on elementary subgroups}

For a Coxeter group $\Ga=\Ga_\Delta$ we define two subschemes: 
$$
\hom_o(\Ga, G)\subset \hom(\Ga, G), \quad \hom_o(\t\Ga, \t G)\subset \hom(\t\Ga, \t G). 
$$
The former consists of homomorphisms   which are injective on every elementary subgroup of $\Ga$; the latter 
 consists of homomorphisms   which are injective on every elementary subgroup of $\t\Ga$ and send $z\in \t\Ga$ to $-1\in SL(2,\C)$. (In fact, the requirement for $z$ follows from faithfulness on elementary subgroups, except when $\Delta$ has no edges.) 
Since elementary subgroups of $\Ga$ and $\t\Ga$ are finite, both $\hom_o(\Ga, G)$ and
$ \hom_o(\t\Ga, \t G)$ are open subschemes of the respective representation schemes. We will see later on that these subschemes are also closed. For each 
$$
\rho\in \hom_o(\t\Ga, G)
$$
we have $\rho(z)=1$, while for each $\t\rho\in  \hom_o(\t\Ga, \t G)$ which projects to $\rho\in \hom_o(\Ga, G)$ we have $\t\rho(z)=-1$. 

\begin{figure}[!ht]
\centering
\includegraphics[width=4in]{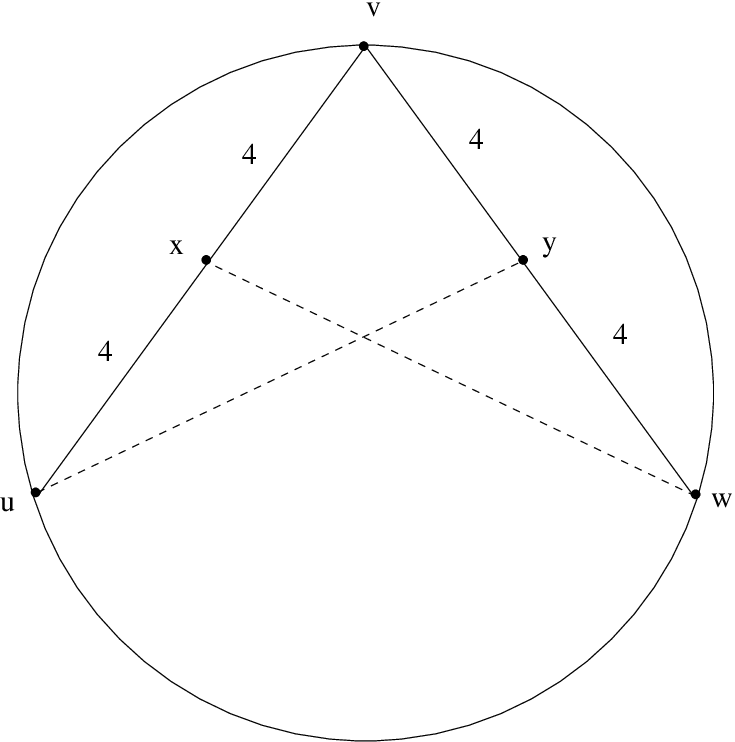}
\caption{Graph $\Omega$.}\label{fig:1}
\end{figure}

In the paper we will be using the labelled graph $\Om$ depicted in the  Figure \ref{fig:1}: This graph has five vertices and nine edges. The edges left unlabeled in the figure, all have the label $2$. The vertices $x, y$ are the only ones not connected to each other by an edge.  

In what follows, we will also use the subgraph $\Upsilon\subset \Om$,  which is the complete graph on the vertices $u, v, w$. The parabolic subgroup $\Ga_\Upsilon < \Ga_\Om$ is isomorphic to $\Z_2^3$. 
Since the group $\Ga_\Upsilon$ is finite, the representation scheme $Hom(\Ga_\Upsilon, G)$ is smooth.

\begin{lem}\label{lem:abel}
Each representation $\rho\in  \hom_o(\Ga, G(\C))$ of the group $\Ga=\Ga_{\Upsilon}$ satisfies:

1. The kernel of $\rho$ is generated by the subgroup
$$
 \<g_u g_v g_w\> \cong \Z_2 
$$
and the image of $\rho$ is isomorphic to $\Z_2\times \Z_2$. 

2. The centralizer of the abelian subgroup $\rho(\Ga)<G(\C)$ in the group $G(\C)$ equals the subgroup $\rho(\Ga)$ itself.  

3.  $\hom_o(\Ga, G(\C))$ is the $G(\C)$-orbit of a singleton $\{\rho_\Upsilon\}$. 
\end{lem}
\proof This lemma is also  elementary:  

1. Consider a homomorphism $\rho\in \hom_o(\Ga, G(\C))$. For each element $\ga\in \Ga$ we let $Fix(\ga)$ denote the fixed-point set of $\rho(\ga)$ in $\C P^1$. The condition that all three involutions $\rho(g_u), \rho(g_v), \rho(g_w)$ are distinct, implies that 
the three fixed-point sets $Fix(g_u), Fix(g_v), Fix(g_w)$ are pairwise disjoint. Commutativity of $\rho(\Ga)$ implies that this group preserves the six-point set
$$
F=Fix(g_u)\cup Fix(g_v)\cup Fix(g_w)\subset \C P^1.  
$$
The element $\rho(g_u)$ fixes $Fix(g_u)$, of course, and defines nontrivial involutions of the other two fixed-point sets 
$$
Fix(g_v)\to Fig(g_v), \quad Fix(g_w)\to Fix(g_w).$$
The same applies to $g_v$ and $g_w$. It follows that 
$$
\rho(g_u g_v)|_F= \rho(g_w)|_F.
$$
Hence, 
$$
\rho(g_u g_v)= \rho(g_w) 
$$
and, thus, 
$$
\< g_u g_v g_w\> < \ker(\rho). 
$$
The equality of these subgroups of $\Ga$ follows from the condition that 
$$
\rho\in \hom_o(\Ga, G(\C)).$$ 
This establishes Part 1 of the lemma. 

2. To prove Part 2, note that every $g\in G$ centralizing $\rho(\Ga)$ has to preserve each set $Fix(g_u)$, $Fix(g_v)$, $Fix(g_w)$. After composing $g$ with elements of $\rho(\Ga)$, we achieve that $g$ fixes the set $Fix(g_u)\cup Fix(g_v)$ pointwise. Therefore, $g\in \rho(\Ga)$. This proves 2. 

3. Part 3 note that, by Part 1, the pull-back morphism
$$
\hom_o(\Z^2, G(\C)) \to \hom_o(\Ga, G(\C)) 
$$
induced by the quotient homomorphism
$$
1\to \< g_u g_v g_w\> \to \Ga\to \Z^2
$$
is surjective. Now, the claim follows from Lemma \ref{lem:Z2timesZ2}. \qed

\begin{lem}\label{lem:Om}
1. $\hom_o(\Ga_\Om, G(\C))$ is a single orbit $G(\C)\cdot \rho_\Om$. 

2. The representation $\rho_\Om$ is infinitesimally rigid: $\h^1(\Ga_\Om, sl(2,\C))=0$.  

3. For each $\rho\in \hom_o(\Ga_\Om, G(\C))$, the adjoint action $Ad\rho$ of $\Ga_\Om$ on the Lie algebra $sl(2,\C)$ has no nonzero fixed vectors. 
\end{lem}
\proof 1. Consider $\rho\in \hom_o(\Ga_\Om,G)$. In view of Lemma \ref{lem:abel}, we can assume that the restriction of 
$\rho$ to  the subgroup $\Ga_\Upsilon$ equals the representation $\rho_\Upsilon$. Consider now the dihedral subgroups 
$$
\<g_u, g_x\>, \quad \<g_x, g_v\>
$$
in the group $\Ga_\Om$. It follows from Lemma \ref{lem:I_2(4)} that there are exactly two extensions $\rho_\pm$ of the representation $\rho_\Upsilon$ to the subgroup $\<g_u, g_x, g_v\>$ (which are faifthful on all elementary parabolic subgroups): For both representations, $\rho_\pm(g_x)$ lies in $SO(3,\R)$, its fixed line in $\R^3$ is contained in the span of the fixed lines of 
$\rho_\Upsilon(g_u),  \rho_\Upsilon(g_v)$. 
This fixed line makes the angle $\frac{\pi}{4}$ with the fixed lines of $\rho_\Upsilon(g_u),  \rho_\Upsilon(g_v)$ and is orthogonal to the fixed line of $\rho(g_w)$.  These representations $\rho_\pm$ are conjugate via the element $\rho(g_v)\in SO(3)$. Therefore, after such conjugation, we fix the value $\rho(g_x)$. We next repeat this argument for the subgroup of $\Ga_\Om$ generated by 
$$
\{g_v, g_y, g_w\}.
$$  
Since conjugation via $\rho(g_w)$ does not alter  $\rho(g_x)$, we obtain the required uniqueness statement.  

2. In what follows we will be using the fact that the adjoint representation of $PSL(2,\C)$ is isomorphic to the complexification $V$ of the standard representation of $SO(3, \R)$ on $\R^3$.  We will also use the notation $V$ and $sl(2,\C)$ for the  representation $Ad \rho$ of the group $\Ga_\Om$ (and its subgroups) in the notation for cocycles and cohomology groups. In particular, for each element 
$a$ of 
$$
\{u, v, w, x, y\}
$$
the fixed-point set of $Ad\rho(g_a)$ is a line in $V$, which we will denote by $V^a$. 
An elementary but useful geometric observation is that 
$$
V^x\subset V^u\oplus V^v
$$
while 
$$
V= V^u\oplus V^v\oplus V^w= V^u\oplus V^x\oplus V^w. 
$$

Consider a cocycle $\xi\in Z^1(\Ga_\Om, V)$. Since $\Ga_\Upsilon$ is finite, $\h^1(\Ga_\Upsilon, V)=0$. Since the restriction of $\xi$ to the subgroup $\Ga_\Upsilon$ is a coboundary, by subtracting off a coboundary from $\xi$, we can assume that $\xi$ vanishes in $\Ga_\Upsilon$. Similarly, there exist 
$\al, \be\in V$ such that 
$$
\xi(h)= \al - Ad\rho(h)\al, \quad \forall h\in\<x, u\>
$$
$$
\xi(h)= \be -Ad\rho(h)\be, \quad \forall h\in\<x, w\>. 
$$
It follows that $\al\in V^u, \be\in V^w$. Moreover, by looking at the value $\xi(x)$, we see that  
$$
\al-\be \in V^x. 
$$
Since the lines $V^u, V^x, V^w$ also span $V$, it follows that $\al=\be=0$. Therefore, $\xi(x)=0$. Similarly, $\xi(y)=0$ and, thus, $\xi=0$ on the entire group $\Ga_\Om$. 

3. This follows from the fact that $\rho(\Ga_\Upsilon)$ has no nonzero fixed vectors in $V=sl(2,\C)$.  
\qed 

\begin{cor}
The scheme $Hom_o(\Ga_\Om,G)$ is smooth. 
\end{cor}

\medskip 
From now on, we will be making the following assumption on the labeled graphs $\Delta$ of Coxeter groups $\Ga$: 

\begin{assumption}\label{ass:1}
1. Every label of the graph $\Delta$ is even. 

2. $\Delta$ contains as a full subgraph the graph $\Om$ above.  
\end{assumption}

\begin{prop}\label{prop:L4}
1. The schemes $\hom_o(\Ga, G)$and $\hom_o(\t\Ga, \t G)$ 
are clopen\footnote{I.e. closed and open.} 
subschemes in $\hom(\Ga, G)$ and $\hom(\t\Ga, \t G)$
respectively. 

2. There is a morphism of schemes $q: \hom_o(\t\Ga, \t G) \to \hom_o(\Ga, G)$, such that for every $\t\rho\in 
\hom_o(\t\Ga, \t G) $ and $\rho=q(\t\rho)$ we have
$$
p\circ \t\rho = \rho \circ \eta, 
$$
where $\eta: \t\Ga\to \Ga$ is the quotient map from \eqref{eq:cen}. 

3. The morphism $q$ is a regular \'etale covering with the deck-group
$\Z_2^{r}$, where $r$ is the rank of $\Ga$.  
\end{prop}
\proof 1. We will give a proof for $\hom_o(\Ga, G)$, since the other  statement is similar. 
Consider an elementary subgroup $\Ga_\Si\subset \Ga$; this subgroup is finite. 
In Lemma \ref{lem:L1} we proved that each irreducible component of 
$\hom(\Ga_\Si, G)$ is a clopen subscheme of $\hom(\Ga_\Si, G)$; furthermore, each component is a single $G$-orbit of a representation $\Ga_\Si\to G$. Then 
$$
\hom_o(\Ga_\Si, G) = \hom(\Ga_\Si, G)\setminus \bigcup_{\theta\in \Ga_\Si - \{1\}} \hom_{\<\theta\>}(\Ga_\Si, G)
$$
is an open subscheme in $\hom(\Ga_\Si, G)$. It is also closed since every subscheme removed was open. 

For each elementary subgroup  $\Ga_\Si<\Ga$ and the inclusion map $\iota_\Si: \Ga_\Si\to \Ga$, 
we have the pull-back morphism
$$
\iota^*_\Si:  \hom(\Ga, G) \to \hom(\Ga_\Si, G). 
$$
Then we have the finite intersection, taken over all elementary subgraphs $\Si\subset \Delta$:   
$$
\hom_o(\Ga, G)= \bigcap_{\Si} (\phi^*_\Si)^{-1}\left(\hom_o(\Ga_\Si, G)\right).
$$
Therefore, $\hom_o(\Ga, G)\subset \hom(\Ga, G)$ is clopen as a finite intersection of clopen subschemes.

\medskip 
2. For each $\t\rho\in \hom_o(\t\Ga, \t G(\C))$, the reduction modulo centers of $\t\Ga$ and $\t G$ yields a 
homomorphism $\rho\in  \hom_o(\Ga, G(\C))$. We need to check that the map 
$$
q: \hom_o(\t\Ga, \t G(\C)) \to \hom_o(\Ga, G(\C)), \quad q(\t\rho)=\rho
$$
obtained in this fashion comes from a morphism of schemes. First, the composition
$$
\t\rho\to p\circ \rho, \quad \hom(\t\Ga, \t G)\to \hom(\t\Ga, G)
$$
is clearly a morphism of schemes. For $\Theta=\{z\}$, we obtain an isomorphism of the schemes 
$$
\hom_\Theta(\t\Ga, G) \to \hom(\Ga, G), 
$$
see Lemma \ref{lem:L3}; and $\hom_\Theta(\t\Ga, G)$ contains the image of 
$\hom_o(\t\Ga, \t G)$. Therefore, $q$ is a composition of two morphisms. 

We next verify surjectivity. Let $\rho\in \hom_o(\Ga, G)$. Define $\t\rho: \t\Ga\to \t G$ by sending generators $g_v$ to arbitrary elements of $p^{-1}(\rho(g_v))$ and sending the central element $z\in \t\Ga$ to $-1\in SL(2, \C)$. In view of Proposition \ref{prop:P0}, for each elementary subgroup $\Ga_\Si$ in $\Ga$, the restriction of $\t\rho$ to the generators of $\t\Ga_\Si$ extends to a faithful homomorphism 
$\t\rho|_{\Ga_\Si}$. 

Since all the relators of the group $\t\Ga$ come from elementary subgroups, it follows that our map of the generators of $\t\Ga$ to $SL(2)$ extends to a homomorphism $\t\rho: \t\Ga\to SL(2)$. This homomorphism belongs to $\hom_o(\t\Ga, \t G)$ since it is faithful on each elementary subgroup. 

Thus, we obtained a surjective morphism 
$$
q: \hom_o(\t\Ga, \t G) \to \hom_o(\Ga, G), q(\t\rho)=\rho. 
$$
The group $\Z_2$  is the group of automorphisms of the covering $\t G\to G$; therefore, the product of $r$ copies of $\Z_2$ 
acts naturally on the product of $r$ copies of $\t G$ as the group of automorphisms of the (regular) cover 
$$
\hat{p}=p\times \ldots \times p: \prod_{i=1}^r \tilde{G}  \to \prod_{i=1}^r {G}. 
$$
Since the rank $r$ of the group $\Ga$ is the number of its generators $g_v$, we have the morphism 
$$
\prod_{i=1}^r \tilde{G}\cong \hom(F_r, \tilde{G}) \stackrel{\hat{p}}{\longrightarrow}  \prod_{i=1}^r {G} \cong \hom(F_r, {G}), 
$$
where $F_r$ is the free group on $r$ generators. We also have the following commutative diagram 
\begin{diagram}
\hom(\t\Ga, \t G) &  \rTo^{q~~~~}  &  \hom(\Ga, G)\\
\dTo &   ~   &  \dTo  \\
\hom(F_r, \tilde{G}) &  \rTo^{\hat{p}~~~~}  &  \hom(F_r, {G})\\
\end{diagram}
where the vertical arrows are the inclusions of representation schemes induced by the epimorphisms
$$
F_r\to \t\Ga, \quad F_r\to \Ga
$$
sending the free generators of $F_r$ to the generators $g_r$ of extended Coxeter and Coxeter groups. 
It is elementary and left to the reader to verify that the group $\Z_2^r$ of automorphisms of 
$\hat{p}$ preserves the subscheme $\hom_o(\t\Ga, \t G)$. Therefore, this finite group acts simply transitively 
on the fibers of the morphism $q$.

3. It remains to show that the map $q$ is etale, i.e., is an isomorphism of analytic germs at every point. 
Let $\t\rho$ be in $\hom_o(\t\Ga, SL(2,\C))$ and set $\rho:=q(\t\rho)$.

Below is a proof which assumes reader's familiarity with \cite{GM}, where the theory of controlling 
differential graded Lie algebras for various deformation problems was developed. 

In view of \cite[Theorem 6.8]{GM}, it suffices to verify that the differential graded Lie algebras controlling these germs are quasi-isomorphic. First, the Lie algebras of $G$ and $\t{G}$ are isomorphic under the covering $p$, which implies that the covering map $p$ induces an isomorphism
$$
\h^i(\t\Ga, Ad\circ \t\rho) \to \h^i(\t\Ga, Ad\circ p(\t\rho)), i\ge 0.  
$$ 
Since the central subgroup $\Z_2$ of $\t\Ga$ is finite, 
$$
\h^i(\Z_2, sl(2,\C))=0, \quad i\ge 1. 
$$
Therefore, applying the Lyndon--Hochshild--Serre spectral sequence to the central 
extension \eqref{eq:cen}, we obtain isomorphisms 
\begin{equation}\label{eq:iso1}
\h^i(\Ga, Ad\circ \rho) \to \h^i(\t\Ga, Ad\circ \t\rho), \rho=q(\t\rho), i\ge 1.  
\end{equation}
(Actually, for $i=0$ both cohomology groups vanish, which implies that they are also isomorphic.)  
These isomorphisms ensure that the morphism 
$$
q: (\hom(\t\Ga, \t G), \t\rho) \to (\hom(\Ga, G), \rho)
$$
is an isomorphisms of germs. \qed 

\begin{rem}
Below is an alternative argument proving that $q$ is etale, which does not reply upon 
differential graded Lie algebras. The morphism $q$ is etale if and only if  $q$ 
induces bijections of sets of $A$-points of representation schemes for all local Artin  $\C$-algebras $A$. 
Let $A$ be a local   Artin  $\C$-algebra and 
$\eps: A\to \C$ be the quotient by the maximal ideal. 
Then we have natural bijections $\hom(\Ga, G(A))\cong \hom(\Ga, G)(A)$ and  
$\hom(\t\Ga, \t{G}(A))\cong \hom(\t\Ga, \t{G})(A)$, and the commutative diagram 
\begin{diagram}
1&\rTo& \hom(\t\Ga,\t{K}_A)    &\rTo & \hom(\t\Ga, \t{G}(A)) & \rTo & \hom(\t\Ga, \t{G}(\C)) & \rTo & 1\\
 &   &  q_{K,A}\dTo  &    &   q_A\dTo       &     &  q\dTo           &     &   \\
1&\rTo& \hom(\Ga,{K}_A)  &\rTo & \hom(\Ga, G(A)) & \rTo & \hom(\Ga, G(\C)) & \rTo & 1\\
\end{diagram} 
where $\t{K}_A $ and $K_A$ are the respective kernels of the group homomorphisms 
$$
\t{G}(A)\to \t{G}(\C), \quad G(A)\to G(\C)
$$
induced by $\eps: A\to \C$.  We observe that the group $\t{K}_A$ is torsion free nilpotent and, since the covering map $p:\t G\to G$ is etale, the induced map $p_{K,A}: \t{K}_A\to K_A$ is an isomorphism. Since $\t{K}_A$ is torsion free, each homomorphism $\t\rho_A: \t\Ga\to  \t{G}(A)$ descends  to a canonical homomorphism $\rho_A: \Ga\to G(A)$. The fact that $\t{K}_A$ is torsion-free also implies that if $\t\rho_A, \t\rho'_A$ have the same projection to $\hom(\t\Ga, \t{G}(\C))$ and define the same homomorphism $\rho_A: \Ga\to G(A)$, then $\t\rho_A= \t\rho'_A$. Therefore, $q_A$ is a injective. The proof that $q_A$ is surjective is also elementary and is left to the reader.    
\end{rem}

\subsection{Character schemes of representations faithful on elementary subgroups}

In this section we extend the results of the previous section 
from representation schemes to character schemes.

\subsubsection{Stability} 

Given a reductive affine algebraic group $H$ and a finitely-generated group $\La$, we have the algebraic action of the group $H$ on the homomorphism scheme $\hom(\La, H)$, 
$$
(h,\rho)\mapsto Inn(h)\circ \rho,
$$
where $Inn(h)$ is the inner automorphism $g\mapsto hgh^{-1}$ of the group $H$. Recall that the character scheme is defined as the Mumford quotient 
$$
X(\La, H)=\hom(\La, H)//H. 
$$
Geometrically speaking, Mumford quotient is obtained by identifying semistable points 
$\hom_{ss}(\La, H)$ of the $H$-action by the {\em extended orbit equivalence relation}, while the restriction of the projection 
$$
\mu: \hom_{ss}(\La, H)\to X(\La, H)
$$
to the stable locus $\hom_{st}(\La, H)$ (consisting of stable points) is just the quotient by the $H$-orbit equivalence. Hence the restriction of the projection to the stable locus has especially simple form. 
We will use the notation
$$
\rho\mapsto [\rho]
$$
for the projection $\mu$. 

A  {\em sufficient condition} for stability of representations   
$\rho\in \hom(\La, H)$ (under the $H$-action) in terms of the Zariski closure of $\rho(\La)$ in $H$
was established in \cite{Johnson-Millson}: 

\begin{thm}
A representation $\rho\in \hom(\La, H)$ is semistable provided that the Zariski closure $\ol{\rho(\La)}$ is  reductive. A representation is stable provided that the Zariski closure $\ol{\rho(\La)}$ is  reductive and the centralizer $Z_H(\rho(\La))$ of the image of $\rho$ is finite.  
\end{thm}

In the case of representations into $H=PO(3)$ and $H=Spin(3)$, the sufficient condition for stability amounts to requiring that the image of $\rho$ is not contained in a Borel subgroup of $H$. Our next goal is to verify stability condition and identify centralizers of the images of representations in the context of Coxeter and extended Coxeter groups. Recall that we are using the notation $G$ for $PSL(2)$ and $\t G$ for $SL(2)$ (regarded as group-schemes). 

\begin{lem}\label{lem:L13}
Let $\Ga$ be a Coxeter group and $\t\Ga$ the corresponding extended Coxeter group,  satisfying 
the Assumption \ref{ass:1}. Then for each $\rho\in \hom_o(\Ga, G)$ and $\t\rho\in \hom_o(\t\Ga, \t G)$ we have: 

1. The representations $\rho, \t\rho$ are stable points in Mumford's sense. 

2. The centralizers of the images of $\rho, \t\rho$ equal the center of the target group.  
\end{lem}
\proof Recall that we require the group $\Ga$ to contain a subgroup $\Ga_\Om$. It suffices to prove both 1 and 2 for the representations 
$\rho\in \hom_o(\Ga_\Om, G), \t\rho\in \hom_o(\t\Ga_\Om, \t G)$, since we have to verify that the image of the representation is not contained in a Borel subgroup and that its centralizer equals the center of the target group. 

(i) First, we consider the case of representations $\t\rho: \t\Ga=\t\Ga_\Om\to SL(2,\C)$. 
We restrict our attention to the subgraph $\Si\subset \Om$, which consists of two vertices $x, y$ and the edge $e=[x,y]$ labelled $4$. Each representation $\t\rho\in \hom_o(\t\Ga, \t G)$ projects to a faithful representation 
$$
\rho: \Ga_\Si\embed PSL(2,\C). 
$$
By Proposition \ref{prop:P0}, the centralizer of the subgroup $\t\rho(\t\Ga_\Si)$ equals the center of $SL(2,\C)$. 
Moreover, the images of the generators of $\t\Ga_\Si$ under $\t\rho$ have distinct eigenlines. It follows that the subgroup $\t\rho(\t\ga_\Si)$ cannot have an invariant line in $\C^2$, thereby proving that $\t\rho(\t\Ga)$ is not contained in a Borel subgroup of $SL(2,\C)$. This proves Parts 1 and 2 for representations to $SL(2,\C)$.  

(ii) Consider now representations $\rho: \Ga\to G$. By the assumption, $\rho$ sends distinct generators of $\Ga$ to distinct elements of $G$. It follows that the group $\rho(\Ga)$ cannot fix a point in $\C P^1$. In other words, the group $\rho(\Ga)$ is not contained in a Borel subgroup of $G$. This proves Part 1.

\medskip 
To prove Part 2, we will use subgroups $\Ga_\Upsilon$ and $\Ga_\Om$ of the group $\Ga$. 
Since $\rho$ belong to $\hom_o(\Ga, PSL(2,\C))$, the centralizer of $\rho(\Ga_\Upsilon)$ in $G$ equals the subgroup $\rho(\Ga_{\Upsilon})$ itself (Lemma \ref{lem:abel}). On the other hand, $\rho$ is faithful on the subgroups generated by 
$\{g_u, g_x\}, \{g_v, g_y\}, \{g_w, g_z\}$. Therefore,  
$$
[\rho(g_u), \rho(g_x)]\ne 1, \quad [\rho(g_v), \rho(g_y)]
\ne 1, \quad [\rho(g_w), \rho(g_z)]\ne 1, 
$$
and, hence, the subgroup $\rho(\Ga)$ has trivial centralizer in $PSL(2,\C)$. \qed

\subsubsection{Cross-sections} 

Let  $Y$ be an (quasiaffine) scheme and $G\acts Y$ be a stable (in Mumford's sense) algebraic group action on $Y$. Suppose that $S\subset Y$ is a closed subscheme, such that the orbit map 
$$
G\times S\to Y
$$
is an isomorphism. (In particular, $S$ projects isomorphically onto $Y//G$). Then $S$ is called a {\em cross-section} for the action of $G$ on $Y$.

\begin{lem}\label{lem:cross-section} 
Suppose that $Y$ is a (quasiaffine) scheme of finite type, $G\times Y\to Y$ is an (algebraic) action of an affine algebraic group, $S\subset X$ is a cross-section for this action. Suppose that $Y'\subset Y$ is a $G$-invariant subscheme. Then $S'=Y'\cap S$ is also a cross-section for the action $G\times Y'\to Y'$. 
\end{lem}
\proof We need to show that the orbit map $G\times S'\to Y'$ is an isomorphism. It suffices to show that for each commutative ring $A$, the orbit map
$$
\mu': G(A)\times S'(A)\to Y'(A)
$$ 
of $A$-points is a bijection. We have $S'(A)=S(A)\cap Y(A)$. Since the orbit map 
$$
\mu: G(A)\times S(A)\to Y(A)
$$  
is a bijection and $Y'(A)$ is $G(A)$-invariant, it follows that  $\mu$ is a bijection. \qed 

\medskip
Note that if the scheme $Y$ and its subscheme $S\subset Y$ are both smooth then the condition that 
$S$ is a cross-section for the action of $G$ is easier to check: It suffices to verify that the set of complex points of $S$ is a set-theoretic cross-section for the action of $G$ and that $G$-orbits are transversal to 
$S$: For each $y\in S(\C)$, 
$$
T_y Y\cong T_y(G y) \oplus T_y S. 
$$

\medskip
We now specialize to the case of representation schemes. Let $\pi'=\pi/N$ be a finitely-generated group (where $\pi$ is a finitely generated group and $N\triangleleft \pi$ is a normal subgroup), 
$G$ is an affine algebraic group $G\times Hom(\pi,G)\to Hom(\pi, G)$ is the action of $G$ by conjugation on the representation scheme. We will think of $Hom(\pi, G)$ as a closed subscheme in the smooth scheme $Hom(F, G)$. Suppose that $U\subset Hom(F,G)$ is a $G$-invariant open affine subscheme and $U'=U\cap Hom(\pi, G)$. We assume that $S\subset U$ is a closed smooth subscheme. Then, in view of smoothness, the property that $S$ is a cross-section for the $G$-action amounts to:

1. $S(\C)$ is a cross-section for the action of $G(\C)$ on $U$. 

2. For each $\rho\in S(\C)$ the action of $\rho(\pi)$ on the Lie algebra of $G(\C)$ has no nonzero invariant vectors. (This condition amounts to the transversality property above.) 

In view of Lemma \ref{lem:cross-section}, the subscheme $S'=S\cap U'$ is a cross-section for the action of $G$ on $U'$. 

\medskip
\subsubsection{Cross-sections of representation schemes} 

We apply the above observations in two situations. First, suppose that $\Ga$ is a Coxeter group satisfying 
the Assumption \ref{ass:1}; we let $G=PO(3)$. We have 
the identity embedding $\iota_\Om: \Ga_\Om\embed \Ga$ of the finite subgroup $\Ga_\Om$. 
Recall that, according to Lemma \ref{lem:Om}, $\hom_o(\Ga_\Om, G(\C))$ consists of a single $G$-orbit 
$G(\C)\cdot \rho_\Om$. We then set
$$
\hom_c(\Ga, G):= (\iota_\Om^*)^{-1}(\rho_\Om). 
$$

 The next lemma is an analogue of Corollary 12.11 in \cite{KM}:

\begin{lem}\label{lem:cross-1}
The subscheme $\hom_c(\Ga, G)$ is  a cross-section for 
the action $G\acts \hom_o(\Ga, G)$. In particular, 
$$
X_o(\Ga, G)\cong \hom_c(\Ga, G). 
$$
\end{lem}
\proof We let $\pi'=\Ga$ and define the new group $\pi$ as the Coxeter group whose Coxeter graph is obtained from the one of $\Ga$ by removing all the edges which are not in $\Om$. Define 
$$
\pi_o:=\Z_2 \star \ldots \star \Z_2. 
$$
Then the representation scheme $Hom_o(\pi, G)$ is smooth, since $\pi$ is isomorphic to the free product 
$$
\pi_o\star \Ga_\Om,
$$
and $U=Hom_o(\Ga_\Om, G)$ is smooth by Lemma \ref{lem:Om}. 
Clearly, $\pi'=\pi/N$ for a normal subgroup $N\triangleleft \pi$. We again have the inclusion 
homomorphism $\iota_\Upsilon: \Ga_\Om \to \pi$; the subscheme
$$
S= \hom_c(\pi, G):= (\iota_\Upsilon^*)^{-1}(\rho_\Om)
$$  
is smooth since it is naturally isomorphic to $Hom(\pi_o)$. The fact that $S$ is a cross-section for the action of $G$ on  $U$ follows immediately from Lemma \ref{lem:Om} and observations following Lemma \ref{lem:cross-section}.  Lastly, note that 
$$
U'= \hom_c(\Ga, G) = U\cap \hom(\Ga, G).
$$
Now, the lemma follows from Lemma \ref{lem:cross-section}.  \qed 

\medskip
The second situation when we apply our description of cross-sections is the one of representations of 
extended Coxeter groups $\tilde\Ga$ (again satisfying the Assumption \ref{ass:1}) to the group $\tilde{G}\cong SL(2)$.  The group $\tilde{G}$ does not act faithfully on $\hom(\tilde\Ga, \tilde G)$; this action factors through the action of the group $G=PO(3)$. 

Earlier, we defined the subscheme  $\hom_o(\t\Ga, \t G)\subset \hom(\t\Ga, \t G)$. Set 
$$
\hom_c(\t\Ga, \t G):= q^{-1}( \hom_c(\Ga, G) ). 
$$

\begin{lemma}\label{lem:cross-2}
$\hom_c(\t\Ga, \t G)$ is a cross-section for the action of $G$ on $\hom_o(\t\Ga, \t G)$. 
\end{lemma}  
 \proof We let $\pi'= \t \Gamma$. Similarly to the proof of  Lemma \ref{lem:cross-1}, we 
 define the extended Coxeter group  $\pi$ by eliminating all the edges which are not in 
  the subgraph $\Om$. Then $\pi'$ is isomorphic to a quotient of $\pi$ and the same proof as in 
  Lemma \ref{lem:cross-1} goes through. \qed

\medskip 
\subsubsection {Character schemes} 

We let $X_o(\Ga, G)$ and $X_o(\t\Ga, \t G)$ denote the projections of 
$\hom_o(\Ga, G)$ and $\hom_o(\t\Ga, \t G)$ to the corresponding character schemes. 

In view of Lemmata  Lemma \ref{lem:cross-1} and Lemma \ref{lem:cross-2}, the projections
$$
\hom_o(\Ga, G)\to X_o(\Ga, G), \quad \hom_o(\t\Ga, \t G)\to X_o(\t\Ga, \t G)
$$
are principal fiber bundles with the structure group $G=PSL(2,\C)$: The center of the group $\t G$ acts trivially on $\hom(\t\Ga, \t G)$. We record this as 

\begin{cor}\label{cor:actions}
There exist natural isomorphisms of germs
$$
(\hom_o(\t\Ga, SL(2)), \t\rho) \cong (\hom_c(\t\Ga, SL(2))\times PSL(2), \t\rho\times 1) \cong 
$$
$$
(X_o(\t\Ga, SL(2))\times PSL(2), [\t\rho]\times 1) 
$$
and
$$
(\hom_o(\Ga, PSL(2)), \rho) \cong \hom_c (\Ga, PSL(2))\times PSL(2), \rho\times 1)
\cong $$
$$
(X_o(\Ga, PSL(2))\times PSL(2), [\rho]\times 1). 
$$
\end{cor}

\subsubsection{Adding a free factor}

Let $F_k$ be the free group on $k$ generators. For an arbitrary finitely-generated group $\La$ and an algebraic group $H$ we have an isomorphism of schemes:  
\begin{equation}\label{eq:isom}
\hom(\La\star F_k, H)\cong \hom(\La, H) \times  \hom(F_k, H) \cong \hom (\La, H)\times H^k. 
\end{equation}
This isomorphism is $H$-equivariant, where the action of $H$ by conjugation on the left side and the diagonal action (by conjugations) on the product space on the right side. We will use these isomorphisms in the two spacial cases: $\La=\Ga, H=G$ and $\La=\t\Ga, H=\t G$, where $G=PSL(2,\C)$, $\t G= SL(2,\C)$ and $\Ga, \t\Ga$ are Coxeter and extended Coxeter groups respectively. Then the isomorphisms \eqref{eq:isom} for these groups allow us to define clopen subschemes 
$$
\hom_o(\Ga\star F_k, G)\subset \hom(\Ga\star F_k, G), \quad \hom_o(\t\Ga\star F_k, \t G)\subset \hom(\t\Ga\star F_k, \t G)
$$
as the images of
$$
\hom_o(\Ga, G)\times G^k , \quad  \hom_o(\t\Ga, \t G)\times \t G^k
$$
respectively. 

It follows from Lemmata \ref{lem:cross-1} and \ref{lem:cross-2} that $\hom_c(\Ga, G) \times G^k$ is a  cross--section for the action of $G$ on $\hom_o(\Ga, G)\times G^k$, while   $\hom_c(\t\Ga, \t G) \times \t{G}^k$ is a  cross--section for the action of $G$ on $\hom_o(\t\Ga, \t G)\times \t{G}^k$. We, thus, obtain:

\begin{lem}\label{lem:L5}
$$
(\hom_o(\Ga, G) \times G^k)/G\cong X_o(\Ga, G) \times G^k$$ 
\end{lem}

The etale covering $q$ defined above yields, for each $k$, the etale covering  
$$
q_k: \hom_o(\t\Ga, \t G)\times \t{G}^k \cong \hom_o(\t \Ga \star F_k, \t G) \to \hom_o(\Ga \star F_k, G) \cong \hom_o(\Ga, G)\times G^k. 
$$

\begin{cor}\label{cor:C6}
1. $X_o(\Ga \star F_k, G)\cong X_o(\Ga, G)\times G^k$. 

2. $\hom_o(\t \Ga \star F_k, \t G)\cong X_o(\t\Ga, \t G) \times \t G^k$. 

3. The covering $q_k$ is \'etale. 
\end{cor}

\section{Universality theorem of Panov and Petrunin}\label{sec:PP}

Proof of Theorem \ref{main} and its corollaries hinges upon two results, the first of which is the following: 

\begin{thm}[Panov--Petrunin Universality Theorem, \cite{PP}]  \label{thm:PP}
Let $\Ga$ be a finitely--pre\-sented group. Then there exists 
a closed 3-dimensional (non-orientable) hyperbolic orbifold $O$ so that $\pi_1(Y)\cong \Ga$, 
where $Y$ is the underlying space of $O$. Furthermore, $Y$ is a 3-dimensional pseudomanifold without boundary. 
\end{thm}

\begin{rem}
Examination of the proof in \cite{PP} shows that the orbifold $O$ admits a hyperbolic manifold cover $\t O\to O$ with the deck-transformation group $\Z_2^4$.  
\end{rem}

The singular set of the pseudomanifold $Y$ consists of singular points $y_j, j=1,\ldots, 2k$, whose neighborhoods $C_j$ in $Y$ are cones over $\R\P^2$. Note that, since $\R\P^2$ has Euler characteristic $1$, the number of conical singularities has to be even. Observe also that one needs $k\ge 1$ in this theorem, since fundamental groups of 3-dimensional manifolds are very restricted among finitely-presented groups. For instance, there are no 3-manifolds $M$ with $\pi_1(M)\cong \Z^4$; therefore, for $\Ga\cong \Z^4$, one cannot have $k= 0$ in Theorem \ref{thm:PP}. 

\begin{problem}
Does Theorem \ref{thm:PP} hold with $k=1$? 
\end{problem}

Given $\Ga$ and $Y$ as in Theorem \ref{thm:PP}, we will construct a closed (non-orientable) 3-dimensional manifold 
$M=M_\Ga$ as follows. (Formally speaking, this 3-manifold also depends on the choice of an orbifold $O$ in Theorem \ref{thm:PP}, which is very far from being unique, however, in order to simplify the notation, we will suppress this dependence).  

Let $O$ be a 3-dimensional orbifold as in Theorem \ref{thm:PP} and let $Y$ be the underlying space of $O$. Let $Y'$ be obtained by removing open cones $C_j, j=1,...,2k$, from $Y$. Then $Y'$ is a compact 3-dimensional manifold with $2k$ boundary components each of which is a copy of the projective plane $\R \P^2$. We let $\theta_i$ denote the generator of the fundamental group of the projective plane $P_i\cong \R\P^2\subset \D M$, which equals the boundary of the cone 
$C_i$. We will regard $\theta_i$ as an element of $\pi_1(Y')$. Set
$$
\Theta:= \{\theta_1,\ldots, \theta_k\}. 
$$
Then
$$
\Ga=\pi_1(Y)=\pi_1(Y')/\<\<\theta_1,\ldots, \theta_{2k}\>\>. 
$$
Next, let $M$ be the closed 3-dimensional manifold obtained by attaching $k$ copies of the product $\R\P^2\times [0,1]$ to $Y'$ along the boundary projective planes, pairing the projective planes $P_i$ and $P_{i+k}$, $i=1,\ldots,k$. 
Then $\pi=\pi_1(M_\Ga)$ is the iterated HNN extension of $\pi_1(Y')$ with stable letters $t_1,...,t_k$:
$$
\left( \left( \left(\pi_1(Y')\star_{\<\theta_1\>}\right) \star_{\<\theta_2\>}\right)...\right)\star_{\<\theta_k\>}. 
$$

Taking the quotient
\begin{equation}\label{eq:phi}
\phi: \pi \to \pi/\<\< \Theta\>\>,
\end{equation}
we, therefore, obtain the group
$\Ga \star F_k$, where $F_k$ is the free group on $k$ generators, projections of the stable letters $t_i, i=1,\ldots, k$ in the above HNN extension. We let
$$
\psi: \Ga \star F_k \to \Ga
$$
denote the further projection to the first direct factor  and set
\begin{equation}\label{eq:xi}
\xi:= \psi\circ \phi: \pi\to \Ga. 
\end{equation}

Now, given an algebraic group $H$, we obtain 
$$
\hom_\Theta(\pi, H)= \phi^*( \hom(\Ga\star F_k, H)),
$$
a clopen subscheme in $\hom(\pi, H)$ (see Corollary \ref{cor:C3}). The isomorphism  
$$
H^k \times \hom(\Ga, H) \stackrel{\cong}{\lra} \hom(\Ga\star F_k, H) \stackrel{\phi^*}{\lra} \hom_\Theta(\pi, H)
$$
restricted to $1\times \hom(\Ga, H)$  equals $\xi^*$. We, thus, we obtain:

\begin{lem}\label{lem:hat}
For each open subscheme $S\subset  \hom(\Ga, H)$, there exists an open subscheme 
$\widehat{S}\subset \hom(\pi, H)$ isomorphic to $H^k\times S$ via the morphism $\phi^*$. 
Furthermore, $\widehat{S}$ contains 
$\xi^*(S)$. 
\end{lem}
\proof Take $\widehat{S}= \phi^*(H^k\times S)$, where we identify $H^k\times \hom(\Ga, H)$ with the representation scheme $\hom(\Ga\star F_k, H)$. \qed

\section{A universality theorem for Coxeter groups}\label{sec:KM}

The second key ingredient we need is the following theorem which is essentially contained in  
\cite{KM}. Before stating the theorem we recall (cf. Lemma \ref{lem:cross-1}) 
that the action $G\acts \hom_o(\Ga, G)$ has a cross--section $\hom_c(\Ga, G)\subset \hom_o(\Ga, G)$, i.e., $\hom_o(\Ga, G)$ is $G$-equivariantly isomorphic to the product $X_o(\Ga, G)\times G$. As always, $G=PO(3)$. 

\begin{thm}[M. Kapovich, J. J. Millson]  
\label{thm:KM}
Let $X$ and $x\in X$ be as in Theorem \ref{main}. Then there exists an open subscheme $X'\subset X$ containing $x$, a finitely-generated Coxeter group $\Ga$ (so that every edge of its graph $\Delta$ has label $2$ or $4$) 
and a representation $\rho_c: \Ga\to PO(3, \R)$ with dense image, so that $X'$ 
is isomorphic to an open subscheme $S'\subset X_o(\Ga, G)$. The representation $\rho_c$ belongs to 
$\hom_o(\Ga, PO(3,\R))$. Furthermore, under this isomorphism, $x$ corresponds to $[\rho_c]$. 
 \end{thm}

\begin{rem}\label{rem:KM}
The fact that the $\hom_o(\Ga, G)\cong X_o(\Ga, G)\times G$ (with $\hom_c(\Ga,G)$ containing $\rho_c$ 
serving as a cross-section for the action $G\acts \hom_o(\Ga, G)$) implies that the preimage $R'_o$ of $S'$ in $\hom_o(\Ga, G)$ is isomorphic to $G\times S'= G\times S'_c$. 
As we saw in \S \ref{sec:Coxeter}, the representation $\rho_c$ lifts to a representation
$$
\tilde\rho_c: \tilde\Ga\to SU(2) 
$$
of the canonical central extension $\tilde{\Ga}$ of $\Ga$. 
\end{rem}

\medskip
Since the universality theorems proven in \cite{KM} are somewhat different from the one stated above, we outline the proof of Theorem \ref{thm:KM}. The main differences are that the results of \cite{KM} are about representations of Shephard and Artin groups rather than Coxeter groups. Furthermore, the representation to $PO(3, \R)$ constructed in \cite{KM} has finite image (which was important for \cite{KM}), although the image group does have trivial centralizer in $PO(3, \C)$.

The arguments below are minor modifications of the ones in \cite{KM}. 

\medskip 
{\bf Step 1 (Scheme-theoretic version of Mn\"ev Universality Theorem).} Without loss of generality, we may assume that the rational point $x$ is the origin $0$ in the affine space 
containing $X$. In \cite{KM} we first construct a {\em based projective arrangement} $A$, so that 
an open subscheme $BR_0(A, \P^2)$ in the space of {\em based projective realizations} $BR(A, \P^2)$, is isomorphic to $X$ as a scheme over $\Q$, and, moreover, the {\em geometrization} isomorphism 
$$
X\stackrel{geo}{\longrightarrow} BR_0(A, \P^2)
$$
sends $x\in X$ to a based realization $\psi_0: A\to \P^2$ whose image is the {\em standard triangle}. Furthermore, the images of the points and lines in $A$ under $\psi_0$ are real.  

\begin{rem}
Subsequently, a proof of this result was also given by Lafforgue in \cite{L}, who was apparently unaware of \cite{KM}. 
\end{rem}

\medskip
{\bf Step 2.} An arrangement $A$ is a certain bipartite graph containing a subgraph $T$ (the ``base'') which is isomorphic to the incidence graph of the ``standard triangle'' (also known as ``standard quadrangle''), see \cite[Figure 7]{KM}. The subgraph $T$ has  
5 vertices $v_{00}, v_x, v_y, c_{10}, v_{01}, v_{11}$ corresponding to the ``points'' of the standard triangle 
and 6 vertices $l_x, l_d, l_y, l_{x1}, l_{y1}, l_\infty$ which correspond to the ``lines'' of the standard triangle.  
In \cite[\S 11]{KM} we further modify the bipartite graph $A$: We make the following identification of vertices:
$$
v_{00}\sim l_{\infty}, \quad v_{x}\sim l_y, \quad v_y \sim l_x.$$ 
We also add to $A$ the edges:
$$
[v_{10}, v_{00}], \quad [v_{01}, v_{00}]. 
$$
We will use the upper case notation $V_{00}=\psi_0(v_{00}), V_{x}=\psi_0(v_x)$, 
etc., to denote vectors in $\C^3$ which project to the images under $\psi_0$ of the point-vertices of $T$. 
The choice of this vectors is not unique, of course; we assume that $V_{00}, V_x, V_y$ form a basis and 
\begin{equation}\label{eq:images}
V_{10}=V_{00}+ V_x, \quad V_{01}=V_{00} + V_y, \quad V_{11}=  V_{00}+ V_x + V_y. 
\end{equation}
This is possible due to the incidencies in $\psi_0(T)$.

However, here, unlike \cite{KM}, {\em we will not add the edge $[v_{00}, v_{11}]$.} (The purpose of this edge in 
\cite{KM} was to ensure that certain representation of a Shephard group is finite.) 
 We let $A'$ denote the resulting graph (no longer bipartite). 
We assign labels to the edges of $A'$ as follows: All edges are labelled $2$ except for the two edges 
$$
[v_{10}, v_{00}], \quad [v_{01}, v_{00}],  
$$
which have the label $4$. We then let $\Ga$ denote the Coxeter group corresponding to this labelled graph. We let $T'$ denote the labelled subgraph of $A'$, whose vertices are the images of the vertices of the arrangement $T$. 

The labeled graph $\Om$ as in Fugure \ref{fig:1} embeds into $T'$ via the map
$$
v\mapsto v_{00}, \quad x\mapsto v_{10}, \quad y\mapsto v_{01}, \quad u\mapsto v_{x}, \quad w\mapsto v_{y}. 
$$

We equip the vector space $\C^3$ with a nondegenerate bilinear form, so that:

1. All subspaces which appear in the image $\psi_0(T)=\psi_0(A)$ are anisotropic (the bilinear form has nondegenerate restriction to these subspaces).

2. The vectors $V_{00}, V_x, V_y\in \C^3$ are pairwise orthogonal and have unit norm. 

\medskip 
We let $PO(3)$ denote the projectivization of the orthogonal group $O(3)$ preserving this bilinear form. 

A realization $\psi\in R(A, \P^2)$ is {\em anisatropic} if for each vertex $v\in A$, the image $\psi(v)$ is an anisatropic subspace in $\C^3$. We will use the notation $R_a(A, \P^2)\subset R(A, \P^2)$ 
and $BR_a(A, \P^2)\subset BR(A, \P^2)$ for open schemes of anisotripic realizations and anisotropic based realizations. By the condition (1) on the inner product above, $BR_a(A, \P^2)$ containins $\psi_0$. 

To every {anisotropic realization} $\psi\in R(A, \P^2)$, we associate a representation of the group $\Ga$ by sending every generator $g_v\in \Ga$ to the isometric involution in $PO(3)$ fixing the subspace $\psi(v)$ in 
$\P^2$. As in \cite{KM}, this map of generators of $\Ga$ to $PO(3)$ defines a representation
$$
\rho_\psi: \Ga\to PO(3,\C).  
$$
We define
$$
\rho_c:=\rho_{\psi_0}. 
$$

By the construction, each representation $\rho_\psi$ is faithful on elementary subgroups: For the edges $[v,w]$ 
in $A$ (where $v$ is a point and $w$ is a line), the incidence condition $\psi(v)\in \psi(v)$ in $P^2$ forces 
the point reflection in $\psi(v)$ be distinct from the line reflection in $\psi(w)$. For the edges 
$$
[v_{10}, v_{00}], \quad [v_{01}, v_{00}] 
$$
the condition \eqref{eq:images} forces the point reflections in $\psi(v_{00}), \psi(v_{10}), \psi(v_{01})$ 
to be pairwise noncommutting and, hence, both subgroups 
$$
\rho_\psi(\<g_{v_{00}}, g_{v_{10}}\>) < PO(3,\C)
$$
$$
\rho_\psi(\<g_{v_{00}}, g_{v_{01}}\>) < PO(3,\C)
$$
are isomorphic to $I_2(4)$. We also note that 
\begin{equation}\label{eq:restriction} 
\rho_\psi|_{\Ga_\Om}= \rho_\Om:= \rho_{\psi_0}: \Ga_\Om\to PO(3,\C). 
\end{equation}
We, thus, obtain the {\em algebraization} morphism of schemes 
$$
alg: BR_a(A, \P^2)\to \hom(\Ga, PO(3)), \quad  \psi\mapsto \rho_\psi. 
$$
As in \cite{KM}, the morphism $alg$ is an isomorphism to its image. It follows from Lemma \ref{lem:cross-1} and \eqref{eq:restriction} that 
the subscheme 
$$
S_c:=alg(BR_a(A, \P^2))\subset \hom_c(\Ga, PO(3))\subset \hom_o(\Ga, PO(3))$$ 
is a cross-section for the action of $G$ on the $G$-orbit of  $S_c$.

\medskip 
Let $\Si\subset A'$ denote the complete subgraph whose vertices are the vertices (points and lines) 
of the standard triangle in $A$, except for the vertex $v_{11}$. As in \cite{KM}, the image under $\rho_c$ of the corresponding parabolic Coxeter subgroup $\Ga_\Si\subset \Ga$, is isomorphic to the finite Coxeter group $B_3$ (the symmetry group of the regular octahedron) divided by the center $\Z_2$. Such a group is a maximal finite subgroup of $PO(3, \R)$. However, the involution $\rho_{c}(g_{v_{11}})$ does not belong to the group 
 $\rho_{c}(\Ga_\Si)$ (this would be order 2 rotation in the center of a face of the octahedron). Thus, the group $\rho_{c}(\Ga)$ has to be dense in $PO(3, \R)$, as it contains (actually, equals to) the dense subgroup
 $ \rho_c(\Ga_{T'})$.   This is the only essential  difference between the construction in this paper and in \cite{KM}, where it was important for the group  $\rho_{c}(\Ga)$ to be finite. 

We let  
$$
\mu: \hom_o(\Ga, G)\to X_o(\Ga, G). 
$$  
denote the restriction of the GIT quotient $\hom(\Ga, G)\to X(\Ga, G)$. Since $\hom_c(\Ga,G)$ is a cross-section 
for the $G$-action on $\hom_o(\Ga, G)$, the morphism $\mu$ is a trivial principal $G$-bundle. 

\begin{thm}\label{thm:based isomorphism} 
$alg: BR_a(A, \P^2)\to \hom_c(\Ga, PO(3,\C))$ is an isomorphism.  
\end{thm} 
\proof We will only sketch the proof since it 
follows closely the argument in Theorem 12.14 in \cite{KM} and the latter is quite long. 
One verifies that $alg$ induces a natural isomorphism of functors  of points. 
For instance, over the complex numbers, each representation $\rho\in \hom_c(\Ga, PO(3,\C))$ 
gives rise to an anisotropic realization:
$\psi(v)\in \P^2(\C)$ is the point fixed by $\rho(g_v)$  (if $v$ is a point-vertex) and 
$\psi(v)\in \P^2(\C)$ is the line fixed by $\rho(g_v)$  (if $v$ is a line-vertex).

\begin{cor}
1. $S_c$ is a cross-section for the action of $G$ on $\hom_o(\Ga, G)$.

2. $\mu\circ alg: BR_a(A, \P^2) \to X_o(\Ga, G)$ is an isomorphism. 

3. $S:=\mu\circ alg(BR_a(A, \P^2))\subset X(\Ga, G)$ is an open subscheme. 
\end{cor}
\proof Part 1 follows from the fact that  
$$
S_c=alg(BR_a(A, \P^2(\C))= \hom_c(\Ga, PO(3,\C))
$$
and the latter is a cross-section for the $G$-action on $\hom_o(\Ga, G)$ (Lemma \ref{lem:cross-1}).  
Part 2 is immediate from Theorem \ref{thm:based isomorphism} and Part 1. Part 3 follows from the fact that $X_o(\Ga, G)$ is an open subscheme in $X(\Ga, G)$. \qed 

\medskip 
We define $X':= geo^{-1}(BR_a(A,\P^2)) \subset X$, an open subscheme in $X$. The composition of $geo, alg$ and $\mu$,  yields an isomorphism $\kappa: X'\to S'$, 
$$
X\supset X' \stackrel{geo}{\lra} BR_a(A,\P^2)\cap BR_0(A, \P^2) 
\stackrel{alg}{\lra} S_c \stackrel{\mu}{\lra} S\subset X_o(\Ga, G),
$$ 
$$
\kappa (X')=S'\subset S,
$$
where $X'\subset X$ and $S'\subset S\subset X_o(\Ga, G)$ are open subsechemes, 
The isomorphism $\kappa$ sends the point 
$x\in X'$ to $[\rho_c]\in X_o(\Ga, G)$. This concludes the proof of Theorem \ref{thm:KM}. \qed 

We let $R_o'$ denote the preimage of $S'$ in $\hom_o(\Ga, G)$ and define  
$$S'_c=  S_c\cap R_o'.$$
Then 
$$
R_o'=G \cdot S'_c \cong G\times S'_c.$$
By inverting the isomorphism $\kappa$ and multiplying with the identity map $G\to G$, we obtain: 

\begin{cor}\label{cor:KM}
There exists an isomorphism of schemes over $\Q$
$$
\omega: R'_o \to G\times X'.  
$$
\end{cor}

\section{Proof of Theorem \ref{main}}

We continue with notation introduced in the previous sections. 
Given an affine scheme $X$ over $\Q$ and a rational point $x\in X$, we use Theorem \ref{thm:KM} to construct a 
Coxeter group $\Ga$ and a representation $\rho_c: \Ga\to PO(3,\R) < PO(3,\C)$. 
Then, as in \S \ref{sec:PP}, we will construct a closed 3-manifold $M=M_\Ga$ with the fundamental group $\pi$, and a clopen subscheme $\hom_o(\pi, G)$ which is isomorphic to the product $\hom_o(\Ga, G)\times G^k$. 
In \eqref{eq:xi} we defined an epimorphism 
$$
\xi: \pi   \stackrel{\phi}{\lra} \Ga \star F_k \stackrel{\psi}{\lra} \Ga. 
$$
Define $\rho_0=\xi^*(\rho_c)\in \hom(\pi, G)$. The subgroup $\rho_0(\pi)= \rho_c(\Ga)< G(\R)$ is dense according to Theorem \ref{thm:KM}.  

We next ``convert'' the open subscheme 
$R'_o\subset \hom_o(\Ga, G)$ to an open subscheme $R'\subset \hom(\pi,G)$. 
The obvious choice $\xi^*(R'_o)$ will not be open in $\hom(\pi,G)$. Instead, we use the open 
subscheme $R':=\widehat{R'_o}\subset \hom(\pi,G)$ defined in Lemma \ref{lem:hat}. 

By combining the isomorphism
$$
R'\to G^k\times R'_o\subset G^k\times \hom(\Ga, G)
$$
with the isomorphism 
$$
id\times \om: G^k \times R'_o\to G^k\times (G \times X') 
$$
(where $\om$ is from Corollary \ref{cor:KM}), we obtain an isomorphism 
$$
f: R' \stackrel{\cong}{\lra} G^k\times R'_o \stackrel{\cong}{\lra} 
 G^{k+1} \times X' \subset G^{k+1}\times X,
$$
sending $\rho_0\in R'$ to 
$$
x'=1 \times x\in G^{k+1} \times X'.$$  
By the construction, $R'$ is open in $\hom(\pi, G)$ and 
$G^{k+1} \times X' $ is open in $G^{k+1}\times X$. 

The cross-section $S'_c\subset Hom_o(\Ga, G)$ in Theorem \ref{thm:KM} yields a cross-section 
$R'_c\subset R\subset Hom(\pi, G)$ for the action $G\acts R'$:
$$
R'_c= \phi^*(G^k \times \psi^*(S'_c)). 
$$
This concludes the proof  
of Theorem \ref{main}. \qed 

\section{Corollaries of Theorem \ref{main}}

Theorem \ref{main} deals with representation schemes of 3-manifold groups 
to $G=PO(3)$; we now consider the corresponding character schemes. 
Since $R'_c\subset Hom_o(\pi, G)$ is a cross-section for the action of $G$ 
on $R'$, Theorem \ref{main} immediately implies:

\begin{cor}
With the notation of Theorem \ref{main}, there exists an open embedding of schemes: 
$$
G^k \times X' \embed X_o(\pi, G)=Hom_o(\pi, G)//G$$
which sends $(1, x)$ to $[\rho_0]$.  In particular, the analytic germ $(\C^{3k}\times X', 0\times x)$
is isomorphic to the analytic germ $(X(\pi, G), [\rho_0])$.    
\end{cor}

We next consider representations of 3-manifold groups to the group $\t G= SL(2)$; we work over $\C$ and, thus, identify $PSL(2,\C)$ with $PO(3,\C)$. Recall that, according to Theorem \ref{thm:KM} 
(and Remark \ref{rem:KM}), for every affine scheme $X$ over $\Q$ and a rational point $x\in X$, 
there exists an open subscheme $X'\subset X$ containing $x$, a Coxeter group $\Ga$ an open subscheme 
$G\times S'_c \cong R'_o\subset \hom(\Ga, G)$,  
and an isomorphism of schemes over $\C$ (which is the identity on the $G$-factor):
$$
R'_o= G\times S'_c\cong G\times S' \to G \times X'
$$
sending $\rho_c\in S'_c\subset \hom_c(\Ga, G) \subset \hom_o(\Ga, G) \subset \hom(\Ga, G)$ to $1\times x$. 

Now, consider representations of the corresponding 
extended Coxeter group $\t \Ga$. Proposition \ref{prop:L4} gives us a $G$-equivariant 
regular  \'etale covering 
$$
q: \hom_o(\t\Ga, \t G) \to \hom_o(\Ga, G)
$$
with the covering group $\Z_2^r$. Restricting to $R'_o\subset \hom_o(\Ga, G)$ we obtain a $G$-equivariant 
regular etale covering
$$
q': \t{R}'_o \to R'_o, \quad \t{R}'_o= q^{-1}(R'_o) \subset   \hom_o(\t\Ga, \t G), \quad q'= q|_{\t{R}'_o}. 
$$
We let $\t \rho_c: \t\Ga\to \t G(\C)$ be a lift of $\rho_c$.  The subscheme $S=\t{R}'_o$ is open 
in $ \hom(\t\Ga, \t G)$. 

\medskip
We now repeat the proof of Theorem \ref{main}. Given the group $\t\Ga$ we construct a closed 3-manifold with the fundamental group $\tilde\pi$ which admits an epimorphism
$$
\tilde\xi: \tilde\pi \to \t\Ga \times F_k \to \t\Ga. 
$$
Set $\tilde\rho_0:= \xi^*(\t\rho_c)$.  

Lemma \ref{lem:hat} (applied to $S$) yields an open subscheme $\widehat{S}\subset \hom(\t\pi, \t G)$ isomorphic to 
$\t G^k \times \t{R}'_o$. Combining this isomorphism with the etale covering
$$
\t G^k \times  \t{R}'_o \to \t G^k \times R'_o
$$
and the isomorphism 
$$
R'_o\cong G \times X' 
$$
we obtain a regular etale covering
$$
\t f: \widehat{S} \to G^{k+1} \times X', 
$$
sending the representation $\t\rho$ to $1\times x$.  The group of covering transformations of $\t f$ is  
the group $\Z_2^r$ (coming from the covering $\t{R}'_o\to R'_o$). 

Furthermore, the subgroup $\t\rho_0(\pi)\subset SL(2, \C)$ is Zariski dense over $\C$, since this subgroup is dense in $SU(2)$ (because it projects to a dense subgroup $\rho_c(\Ga)\subset PO(3, \R)$). Since all groups $p(\t \rho(\pi))$, $\t \rho\in Hom_o(\pi, \t G)$, contain a conjugate of the group 
$$
\rho_c(\Ga_{T'}), 
$$ 
it also follows that for every $\t\rho\in \widehat S$, the group $\t\rho(\pi)$ is Zariski dense in 
$SL(2,\C)$. (The subgraph $T'\subset A'$ is defined in \S \ref{sec:KM}); the group $\Ga_{T'}\subset \Ga_{A'}=\Ga$ is the corresponding parabolic subgroup.) 
This proves Corollary \ref{cor:main}. 
\qed

\section{Orbifold-group representations}

Let $\hat\Ga$ be the fundamental group of the hyperbolic orbifold appearing in Theorem \ref{thm:PP}. 
This group contains cyclic subgroups $\<\theta_i\>\cong \Z_2, i=1,\ldots, 2k$, corresponding to the singular points $y_i$. 
The group $\Ga$ is the quotient
$$
\hat\Ga/\<\<\hat\Theta\>\>,
$$
where $\hat\Theta=\{\theta_1,\ldots, \theta_{2k}\}\subset \hat\Ga$. Then for every algebraic group $H$, 
$$
\hom(\Ga, H)\cong \hom_{\hat\Theta}(\hat\Ga, H) 
$$
and the latter is an open subscheme in $\hom(\hat\Ga, H)$ (see Corollary \ref{cor:C3}). Now, let $\Ga$ be a Coxeter group (as in Theorem \ref{thm:KM}) or its canonical central extension. In view of Theorems \ref{thm:PP} and \ref{thm:KM}, one obtains:

\begin{cor}
Theorems \ref{main} and Corollaries \ref{cor:main}, \ref{cor:germs}, 
also hold for groups $\pi$ which are fundamental groups 
of 3-dimensional closed hyperbolic orbifolds. 
\end{cor}

By passing to a finite-index torsion-free subgroups of $\pi$, in view of \cite[Theorem 5.1]{KM0}, we obtain new examples of fundamental groups of hyperbolic 3-manifolds and their representations to $SO(3), SU(2)$ with non-quadratic singularities of character varieties, cf. \cite[Theorem 5.1]{KM0} where it is proven that nonquadratic singularities of character schemes are inherited by finite index subgroups. (The first such examples were constructed in \cite{KM0}.) 

\begin{question}
Do Theorems \ref{main} and Corollaries \ref{cor:main}, \ref{cor:germs}, 
also hold for groups $\pi$ which are fundamental groups 
of 3-dimensional closed hyperbolic manifolds? Do they hold for 3-dimensional manifolds which are 3-dimensional (integer or rational) homology spheres?  
\end{question}

\medskip
{\bf Older examples.} We note that the first example of a nonreduced representation scheme was constructed by Lubotzky and Magid in \cite[p. 43]{LM}: They start with the von Dyck group  
$$
\Ga= \< a, b, c| a^3=b^3=c^3=abc=1\>
$$
and its representation 
$$
\rho: \Ga\to SL(2, \C)
$$
whose image is a cyclic group $\Z_3$ of order $3$ ($\rho$ sends $a, b$ and $c$ to the same generator of 
$\Z_3$). Then $\h^1(\Ga, Ad\rho)\cong \C$, while 
$[\rho]\in X^{red}(\Ga, SL(2, \C))$ is an isolated point. It follows that $B^1(\Ga, Ad\rho)\cong \C^2$ and, 
therefore,
$$
Z^1(\Ga, Ad\rho)\cong \C^3. 
$$
The same, of course, holds for all representations $\Ga\to SL(2,\C)$ conjugate to $\rho$. 
On the other hand, the component $C$ of $\rho$ in $Hom^{red}(\Ga, SL(2,\C))$ is a smooth 2-fold, isomorphic to $SL(2,\C)/T$, where $T\cong \C^*$ is a maximal torus in $SL(2,\C)$. In contrast, 
the Zariski tangent space $T_c C$ at each $c\in C$ is isomorphic to $Z^1(\Ga, Ad\rho)\cong \C^3$. 
It follows that the representation scheme $\hom(\Ga, SL(2))$ is nonreduced.

This example is promoted to a nonreduced representation scheme of a 3-manifold group as follows. Consider a closed 3-dimensional Seifert manifold $M$ which is an oriented Seifert-bundle over the orbifold $S^2(3,3,3)$, i.e., over the sphere with 3 cone points of order $3$. The fundamental group of $S^2(3,3,3)$ is the von Dyck group $\Ga$ above, while 
$\pi=\pi_1(M)$ is the central extension of $\Ga$:
$$
1\to \Z \to \pi \to \Ga \to 1, 
$$
$$
\pi= \< a, b, c, z| a^3=b^3=c^3= abc= z, [a,z]=[b,z]=[c,z]=1\>. 
$$
The representation $\rho$ lifts to a representation $\t\rho: \pi\to SL(2,\C)$ whose kernel contains the center $\<z\>$ of $\pi$. Then, a direct computation shows that   
$$
\h^1(\pi, Ad \rho)\cong \C^2, \quad Z^1(\pi, Ad \rho)\cong \C^4,
$$
while the germ of the reduced representation variety $\hom^{red}(\pi, SL(2,\C))$ at $\t\rho$ is a smooth 
3-fold, consisting of reducible representations $\pi\to SL(2,\C)$. Up to conjugation, these representations all have the form
$$
\rho_t(a)=\rho_t(b)=\rho_t(c)=
\left[\begin{array}{cc}
t&0\\
0&t^{-1}\\
\end{array}
\right], \quad t\in \C^*. 
$$  

\medskip 
The advantage of the examples constructed in Theorem \ref{main} and its corollaries, is that the representation and character schemes constructed there are nonreduced at points corresponding to representations with trivial centralizer.

\medskip
Addresses:

\medskip 

\noindent M.K.: University of California, Davis, CA 95616 

{\begin{verbatim}kapovich@math.ucdavis.edu\end{verbatim}}

\noindent J.M.: University of Maryland, College Park, MD 20742 

{\begin{verbatim}jjm@math.umd.edu\end{verbatim}}

\end{document}